\title{Economical toric spines via Cheeger's Inequality}
\author{
Noga Alon\thanks{Schools of
Mathematics and Computer Science, Raymond and Beverly Sackler
Faculty of Exact Sciences, Tel Aviv University, Tel Aviv, 69978,
Israel and IAS, Princeton, NJ 08540, USA. Email: nogaa@tau.ac.il.
Research supported in part by a USA-Israeli BSF grant, by the Israel
Science Foundation and by the Hermann Minkowski Minerva Center for
Geometry at Tel Aviv University.}
\and Bo'az Klartag\thanks{Department of Mathematics, Princeton University,
Princeton, NJ 08544, USA and 
School of
Mathematics, Raymond and Beverly Sackler
Faculty of Exact Sciences, Tel Aviv University, Tel Aviv 69978,
Israel. Email: bklartag@princeton.edu.
Research supported in part by a Clay Research Fellowship and by an
NSF grant.}
}
\date{}
\documentclass[11pt]{article}
\usepackage{amssymb}
\usepackage{amsmath}
\usepackage{amsfonts}

\newtheorem{theo}{Theorem}

\newtheorem{lemma}[theo]{Lemma}

\def \RR {\mathbb R}
\def \TT {\mathbb T}
\def \ZZ {\mathbb Z}
\def \EE {\mathbb E}
\def \bP {\mathbb P}
\def \vphi {\varphi}
\def \eps {\varepsilon}

\begin{document}
\maketitle
\begin{abstract}
Let $G_{\infty}=(C_m^d)_{\infty}$ denote the graph whose set of vertices
is $\{1,\ldots ,m\}^d$, where two distinct vertices are adjacent iff
they are either equal or adjacent in $C_m$ in each coordinate. Let
$G_{1}=(C_m^d)_1$ denote the graph on the same set of vertices in which
two vertices are adjacent iff they are adjacent in one coordinate in $C_m$
and equal in all others. Both  graphs  can be viewed as graphs of the
$d$-dimensional torus. We prove that one can delete $O(\sqrt d m^{d-1})$
vertices of $G_1$ so that no topologically nontrivial cycles  remain. This
improves an $O(d^{\log_2 (3/2)}m^{d-1})$ estimate of Bollob\'as, Kindler,
Leader and O'Donnell. We also give a short proof of a result implicit in
a recent paper of Raz: one can delete an $O(\sqrt d/m)$ fraction of the
edges of $G_{\infty}$ so that no topologically nontrivial cycles remain in
this graph. Our technique also yields a short proof of a recent result of
Kindler, O'Donnell, Rao and Wigderson; there is a subset of the continuous
$d$-dimensional torus of surface area $O(\sqrt d)$ that intersects
all nontrivial cycles. All proofs are based on the same general idea:
the consideration of random shifts of a body with small boundary and no-
nontrivial cycles, whose existence is proved by applying the isoperimetric
inequality of Cheeger or its vertex or edge discrete analogues.
\end{abstract}
\section{Introduction}
Let $G_{\infty}=(C_m^d)_{\infty}$ denote the $d$-(AND)-power of the cycle $C_m$ on the
vertices $M=\{1,2, \ldots ,m\}$, 
that is, the graph whose set of vertices is $M^d$, where two distinct
vertices $(i_1, i_2, \ldots, i_d)$ and $(j_1, j_2, \ldots, j_d)$
are adjacent iff for every index $s$, $i_s$ and $j_s$ 
are either equal or 
adjacent in $C_m$. Similarly, let $G_{1}=(C_m^d)_1$ denote the graph
on the set of vertices $M^d$ in which two vertices $(i_1, i_2, \ldots, i_d)$ and $(j_1,
j_2, \ldots, j_d)$ are adjacent iff they are equal in all coordinates but one, in
which they are
adjacent in $C_m$.

Both graphs $G_{\infty}$ and $G_1$ can be 
viewed as graphs of the $d$-dimensional torus.
A cycle in any of them is called
{\em nontrivial} if it wraps around the torus, that is, 
if its projection
along at least one of the coordinates contains the full cycle $C_m$.
A {\em spine} (or an {\em edge-spine})
is a set of edges that intersects
every nontrivial cycle. 
It is easy to see that there is a spine in $G_{\infty}$ containing
a fraction of $O(d/m)$ of the  edges. A recent result of Raz \cite{Ra}, motivated by the 
investigation of parallel 
repetition of the odd cycle game, can be used to show 
that there are much smaller spines consisting of only a fraction of $O(\sqrt d/m)$ of
the edges. Here we prove the following
sharper version of this result.
\begin{theo}
\label{t1}
There exists  an edge-spine of $G_{\infty}$ containing a fraction
of at most 
$2 \mu/(3^d-1)=O(\sqrt d/m)$ of the edges of $G_{\infty}$, where here
$\mu= \sqrt {2 \cdot (3^d-1) \cdot (3^d-(1+2\cos(\pi/m))^d)}.$
\end{theo}
It is not difficult to see that the size of the smallest edge-spine in
$G_1$ is precisely $dm^{d-1}$. Indeed, the set 
$$
\{ ~~\{(i_1, \ldots,i_{s-1},0,i_{s+1}, \ldots  i_d),(j_1, \ldots, j_{s-1},1,
j_{s+1}, \ldots, j_d)\}: ~ 1 \leq s \leq d, i_r, j_t \in C_m~~\}
$$
forms a spine, and there is no smaller spine as the set of all edges of 
$G_1$ can be partitioned into $dm^{d-1}$ pairwise edge disjoint nontrivial cycles .

A {\em vertex-spine} is a set of vertices that intersects every nontrivial cycle.
For vertex spines, the smallest size is known for $G_{\infty}$ and is not known
for $G_1$. Indeed, improving a result of \cite{SSZ}, it is  proved in \cite{BKLO}
that the size of the
smallest vertex spine in $G_{\infty}$ is  $m^d-(m-1)^d$, that is, the vertex-spine consisting
of all vertices in which at least one coordinate is $0$ is of minimum size. 
For $G_1$ the situation is more complicated. 
It is easy to see that there is a vertex spine  consisting of at most $dm^{d-1}$
vertices. This has been improved in \cite{BKLO}, where it is shown that there is a
vertex spine of size at most $d^{\log_2 (3/2)} m^{d-1}\approx d^{0.6} m^{d-1}$.
The following result improves this estimate.
\begin{theo}
\label{t11}
There exists a vertex spine of $G_1$ containing at most $2 \pi \sqrt d m^{d-1}$ vertices.
\end{theo}

The discrete results above have a continuous analogue studied in
\cite{KORW}. 
Let $\TT^d = \left. \RR^d \right / \ZZ^d$ 
be the $d$-dimensional unit torus.
We write $\mbox{Vol}_d$ and $\mbox{Vol}_{d-1}$ for the
$d$-dimensional and $(d-1)$-dimensional Hausdorff measures 
on the unit torus $\TT^d$. 
A loop is a continuous image of the circle. A loop in $\TT^d$ 
is called contractible if it may be continuously deformed  
to a single point
in $\TT^d$.
A {\it spine} in $\TT^d$ is a subset $S \subset \TT^d$ 
that intersects any non-contractible 
loop. Clearly,  the set
$$ S = \{ (x_1,\ldots,x_d) \in [0,1)^d ; 
\exists i \in \{1,\ldots,d \} , x_i = 0 \} $$
is a spine, with $\mbox{Vol}_{d-1}(S) = d$. 
In \cite{KORW} it is shown that we can find  a much smaller spine.
\begin{theo} [\cite{KORW}]
\label{thm_1053}
There exists a compact spine 
$S \subset \TT^d$ with $\mbox{Vol}_{d-1}(S) \leq 2 \pi \sqrt{d}$.
\end{theo}
In this paper we give relatively short proofs of the above three theorems.
The crucial observation is that  in all three cases one can apply either
the isoperimetric inequality of Cheeger (see, e.g., \cite{SY} for a short proof), 
or its discrete version for vertex boundary  (proved in \cite{Al}) or for edge boundary 
(see, e.g., \cite{ASS}), to obtain 
a substructure (an induced subgraph in the discrete case, and a body
in the continuous case) containing no nontrivial cycles,
whose boundary is small with respect to its
volume. The required spine is constructed  in all cases by pieces
of the boundaries of random
shifts of this substructure. The proofs given here, while 
related to the  ones given in \cite{Ra} and \cite{KORW}, are 
significantly shorter.
More importantly, they supply a clear explanation for the choice
of the functions  whose level sets provide  the required 
substructures, as these appear naturally as 
eigenfunctions of the corresponding Laplace operators.
%\iffalse
Indeed, the proof in \cite{KORW} also produces a spine by combining pieces 
of boundaries  of random shifts of  (several) level sets  of an appropriate function, but
provides no clear intuition to the choice of this function. It also estimates the
boundary using a slightly more complicated argument than the one given here,
and thus requires a somewhat
tedious computation. The proof in \cite{Ra} uses a different function, yielding a slightly weaker conclusion.
%\fi
As we briefly remark, our
approach can be applied to derive similar results for
other examples of graphs and bodies.

\section{The discrete case}

\subsection{Edge spines}
We start with a review of the discrete version of Cheeger's inequality
with Dirichlet boundary conditions. For completeness, we include its proof.
\begin{theo}
\label{t3}
Let $G=(V,E)$ be a graph, where $V=\{1,\ldots, n\}$, 
let $A=(a_{ij})_{i,j \in  V}$ be its 
adjacency matrix and let $Q=\mbox{diag}(d(i)_{i \in V})-A$ 
be its Laplace matrix, where $d(i)$ is the degree of $i$. 
Let $U \subset V$ be a set of vertices, and let
$x=(x_1,x_2, \ldots ,x_n)$  be a vector assigning a value $x_i$ to 
vertex number $i$. Assume, further that $x_j=0$ for all $j \in U$ 
and that for every $W \subset V-U$, $e(W,V-W) \geq c |W|$, where
$e(W,V-W)$ is the number of edges joining a vertex of $W$ with one of
its complement. Then $x^t Q x \geq \frac{c^2}{2D} \sum_{i \in V} x^2_i,$
where $D$ is the maximum degree of a vertex of $G$.
\end{theo}
{\bf Proof: }\,
Without loss of generality assume that $V-U=\{1,2,\ldots,r\}$
and that $x_1^2 \geq x_2^2 \ldots \geq x_r^2$. Since
$x^t Q x =\sum_{ij \in E} (x_i-x_j)^2$, by Cauchy Schwartz:
\begin{equation}
\label{e11}
2D (\sum_{i\in V} x_i^2) \cdot x^t Q x  \geq 
\sum_{ij \in E} (x_i +x_j)^2 \sum_{ij \in E} (x_i-x_j)^2
\geq [\sum_{ij \in E, i<j} (x_i^2-x_j^2)]^2.
\end{equation}
Replacing each term $x_i^2-x_j^2$ in the last expression
by $(x_i^2-x_{i+1}^2)+(x_{i+1}^2-x_{i+2}^2)+ \cdots + (x_{j-1}^2-x_j^2)$,
the expression obtained from the sum
$S=\sum_{ij \in E, i<j} (x_i^2-x_j^2)$ contains each term of the form
$x_i^2-x_{i+1}^2$ exactly $e(\{1,2,\ldots ,i\}, \{i+1, \ldots ,n\})$
times, and by assumption this number is at least $ci$ for all $ i \leq r$.
As $x_i=0$ for $i >r$ this implies that $S\geq \sum_{i \leq n} 
ci (x_i^2-x_{i+1}^2)
=c \sum_{i \in V} x_i^2$ (where, by definition, $x_{n+1}=0$).
Plugging in (\ref{e11}) the desired result follows.
$\Box$

\noindent
{\bf Remark:}\, Note that the proof works even if we only assume that
$e(W,V-W) \geq c|W|$ for every $W$ which is a level set of the vector 
$(x_i^2)_{i \in V}$, that is, for every $W$ consisting of all vertices $i$
with $x_i^2 \geq t$. 
\vspace{0.4cm}

\noindent
Let $C_m$ be, as before, the cycle of length $m$ on the set of 
vertices $M=\{1,2, \ldots m\}$ (in this order),
and let $G_{\infty}=(C_m^d)_{\infty}=(V,E)$ denote its $d$-(AND)-power.
Note that $G_{\infty}$
is $D=3^d-1$ regular. 
\begin{lemma}
\label{l4}
There exists a set $W$ of vertices of $G_{\infty}$ that 
contains no nontrivial cycles such that
$e(W,V-W) \leq\mu |W|$, where $\mu$ is as in Theorem \ref{t1},
and
%$\mu= \sqrt {2 \cdot (3^d-1) \cdot (3^d-(1+2\cos(\pi/m))^d)}$
satisfies $\mu / (3^d - 1) = O(\sqrt{d} / m)$.
\end{lemma}

\noindent
{\bf Proof:}\,  Let $A=(a_{ij})_{i,j \in M}$ be the adjacency matrix of
$C_m$, and let $A'$ be the matrix obtained from it by 
replacing the last row and last column  by the zero vector. Note that the
adjacency matrix of $G_{\infty}$ is $(I+A)^{\otimes d} -I^{\otimes d}$, where
for every matrix $B$, 
$B^{\otimes d}$ denotes the tensor product of $d$ copies of $B$. Note also that
if $x$ is a vector of length $m$ and $x_m=0$, and if $x^{\otimes d}$ is the 
tensor product of $d$ copies of $x$, then 
$$
(x^{\otimes d})^t [(I+A)^{\otimes d} -I^{\otimes d}] x^{\otimes d}
=(x^{\otimes d})^t [(I+A')^{\otimes d} -I^{\otimes d}] x^{\otimes d},
$$ 
since these two matrices differ only in entries where the contribution
to the quadratic form  vanishes, as $x_m=0$.

A simple computation shows that the vector $x=\sin(\pi j/m)_{j \in M}$
(which satisfies $x_m=0$) is an eigenvector of $A'$ with eigenvalue
$\lambda=2\cos(\pi/m)$. Therefore, $x^{\otimes d}$ is an eigenvector
of $(I+A')^{\otimes d} -I^{\otimes d}$ with eigenvalue
$\Lambda=(1+\lambda)^d-1=(1+2\cos(\pi/m))^d-1$.  By the above discussion
this implies that
$$
(x^{\otimes d})^t [(I+A)^{\otimes d} -I^{\otimes d}] x^{\otimes d}
=\Lambda || x^{\otimes d} ||^2,
$$ 
and as the Laplace matrix of $G_{\infty}$ is
$Q=(3^d-1) I^{\otimes d}- 
[(I+A)^{\otimes d} -I^{\otimes d}]$ this implies
that
$$
\frac{(x^{\otimes d})^t Q x^{\otimes d}}{||x^{\otimes d}||^2 }
=3^d-1-\Lambda.
$$
By Theorem \ref{t3} we conclude that 
there is a subset $W$ of the vertices
of $G_{\infty}$ that contains no vertex  
with any coordinate being $m$, so that
$ e(W,V-W) \leq \mu |W| $ with $\mu$ as in Theorem \ref{t1}.
The induced subgraph on $W$ contains no nontrivial cycle since
$W \subset (M-\{m\})^d$.    $\Box$
\vspace{0.4cm}

\noindent
{\bf Proof of Theorem \ref{t1}:}\, 
Let $W$ be as in Lemma \ref{l4}, 
let $v_1, v_2, \ldots$ be a random sequence
of vectors in $Z_m^d$, and define $W_i=v_i+W=\{v_i+w: w \in W\}$ where
addition is taken modulo $m$ in each coordinate. 
By symmetry, the induced subgraph of
$G_{\infty}$ on each $W_i$ is isomorphic to that on $W$ and hence contains 
no nontrivial cycle. Obviously, with probability $1$ there exists a finite
$s$ so that $\cup_{i=1}^s W_i=V$. For each $ i $, let $E_i$
be the set of all edges that connect a vertex of $W_i-\cup_{j<i}W_j$ to
a vertex outside $W_i$. The union of all these sets $E_i$ is 
clearly a spine, as each cycle that uses no edge of this union is contained in 
a single set $W_i$. We claim that the expected value of the random variable
$\sum |E_i|$  is at most $\mu m^d$. 
To see this, observe  that by Lemma \ref{l4}, if we choose a random vertex of
$W$ and a random edge incident with it, the probability that  this edge leads to
$V-W$ is at most $\frac{\mu |W|}{(3^d-1)|W|}=\frac{\mu}{3^d-1}$. 
Fix a vertex $v \in M^d$, and 
let $i$ be the smallest $j$ so that $v \in W_j$. Conditioning on
$i$ being the smallest such $j$, $v$ is a uniform random vertex of $W_i$,
and hence if we now choose a random edge incident with it, the probability
it leads to a vertex outside $W_i$ is at most $\frac{\mu}{3^d-1}$.
It follows that the expected size of $E_i$ is at most the expected size
of $W_i-\cup_{j<i}W_j$ times $\mu$. Summing over all values of $i$
and using the fact that with probability $1$ the union of all sets $W_i$ is $V$
we conclude that the expected value of $\sum|E_i|$ is at most 
$\mu |V|=\mu m^d$. Thus there is a choice of sets $W_i$ so that the spine
$\cup_i E_i$ they provide is of size at most $\mu m^d$, 
completing the proof.  $\Box$
%\vspace{0.4cm}
%\newpage
\subsection{Vertex spines}
We need the following version of the inequality of \cite{Al} with Dirichlet boundary
condition. This is an analog of Theorem \ref{t3}, dealing with vertex boundary instead
of edge boundary. Its proof, which is based on the arguments in \cite{Al}, is somewhat more
complicated than that of Theorem \ref{t3}. 
\begin{theo}
\label{t13}
Let $G=(V,E)$ be a graph, where $V=\{1,\ldots, n\}$, 
let $A=(a_{ij})_{i,j \in  V}$ be its 
adjacency matrix and let $Q=\mbox{diag}(d(i)_{i \in V})-A$ 
be its Laplace matrix, where $d(i)$ is the degree of $i$. 
Let $U \subset V$ be a set of vertices, and let
$x=(x_1,x_2, \ldots ,x_n)$  be a vector assigning a value $x_i$ to 
vertex number $i$. Assume, further that $x_j=0$ for all $j \in U$ 
and that for every $W \subset V-U$, $|N(W)-W| \geq c |W|$, where
$N(W)$ is the  set of all vertices that have a neighbor in $W$. 
Then $x^t Q x \geq \frac{c^2}{4+2c^2} \sum_{i \in V} x^2_i.$
\end{theo}
{\bf Proof: }\,
Put $Y=V-U$. We claim that there is an orientation $\overline{E}$ of $E$ and
a function $h: \overline{E} \mapsto [0,1]$ so that the sum
$\sum_{j,(i,j) \in \overline{E}} h(i,j) $ 
is at most $1+c$ for all $i \in Y$,
%and is $0$ for all $i \in U$, 
the sum 
$\sum_{j,(j,i) \in \overline{E}} h(j,i)$ is at most $1$ for all $i$
and the difference
$\sum_{j,(i,j) \in \overline{E}} h(i,j)-\sum_{j,(j,i) \in \overline{E}} h(j,i)$
is at least $c$ for every $i \in Y$.
To prove this claim, consider the network flow problem in which the 
set of vertices consists of a source $s$, a sink $t$, a set $Y'$ consisting
of a copy $y'$ of every $y \in Y$ and a set $V"$ consisting of a copy $v"$ of every
vertex $v \in V$. 
For each $y' \in Y'$, 
$(s,y')$ is an arc of the network with capacity $1+c$, for each
$v" \in V"$, $(v",t)$ is an arc of the network  with capacity $1$, and 
in addition, $(u',u")$ is an arc of capacity $1$ for each $u \in Y$, and
for  each edge $uv$ of $G$, with $u \in Y$, $v \in V$, the arc $(u',v")$
belongs to the
network, and has capacity $1$. (Note that if $v$ is also in $Y$, then the arc
$(v',u")$ is also in the network.)
It is not difficult to check that the value of the maximum
flow in this network is $(1+c)|Y|.$ Indeed, suppose we are given a cut and let $X \subseteq Y$ be the set of all vertices $y \in Y$ such that $(s,y^{\prime})$ belongs to the cut. Then the cut must contain, for each $v \in (Y-X) \cup N(Y-X)$,
at least one arc incident with $v"$. As these arcs are pairwise distinct and there
are at least $(1+c)|Y-X|$ of them, each having capacity $1$, it follows that
the total capacity of the cut  is at least 
$(1+c)|X| +(1+c) |Y-X|=(1+c) |Y|.$ 
By the Maxflow-Mincut Theorem there exists a flow of value  at least 
$(1+c)|Y|$, and this is clearly a maximum flow that saturates all edges
$(s,y')$ with $y \in Y$. 
If there is a positive flow in two arcs $(i',j")$ and 
$(j',i")$ (for some $i,j \in Y$), subtract the minimum of these two from both,
to ensure that at least one of these two quantities is zero, and subtract this minimum
from the value of the flow on $(s,i'),(j",t)$ and on $(s,j'), (i",t)$, thus keeping it
a valid flow without changing the value of the difference between the total flow
leaving $i'$ and the total flow going into $i"$.
Let $h'$ be the resulting flow. If $h'(i',j")>0$ for $ij \in E$,
orient the edge $ij$ from $i$ to $j$ (in case $h'(i',j")=h'(j',i")=0$ orient
the edge arbitrarily). Finally, for each oriented edge $(i,j)$, define 
$h(i,j)=h'(i',j")$.  One can easily check that the function $h$ 
satisfies the assertion of the claim.

We next note that the properties of $h$ imply the following two inequalities.
\begin{equation}
\label{e131}
\sum_{(i,j) \in \overline{E}} h^2(i,j) (x_i +x_j)^2 \leq (4 +2c^2) \sum_i x_i^2.
\end{equation}
\begin{equation}
\label{e132}
\sum_{(i,j) \in \overline{E}} h(i,j) (x_i^2 -x_j^2) \geq c \sum_i x_i^2.
\end{equation}
Indeed, (\ref{e131}) follows, as
$$
\sum_{(i,j) \in \overline{E}} h^2(i,j) (x_i +x_j)^2 
\leq 2 \sum_{(i,j) \in \overline{E}} h^2(i,j) (x_i^2 +x_j^2)
$$
$$
= 2\sum_{i \in Y} x_i^2 (\sum_{j,(i,j) \in \overline{E}} h^2 (i,j)
+ \sum_{j,(j,i) \in \overline{E}} h^2 (j,i))
\leq 2 (2+c^2) \sum_{i} x_i^2,
$$
where here we used the fact that $x_i=0$ for all $i \not \in Y$ and the fact
that the sum of squares of reals in $[0,1]$ whose sum  is at most $(1+c)$ does not exceed
$1+c^2$ (and the sum of squares of real numbers in $[0,1]$ whose sum  is at most $1$ does
not exceed $1$). 

To prove (\ref{e132}) note that
$$
\sum_{(i,j) \in \overline{E}} h(i,j) (x_i^2 -x_j^2)
=\sum_{i \in Y} x_i^2 ( \sum_{j,(i,j) \in \overline{E}} h (i,j)
-\sum_{j,(j,i) \in \overline{E}} h (j,i)) \geq c \sum_i x_i^2.
$$

We can now complete the proof of the theorem. Indeed, by Cauchy-Schwartz,
(\ref{e131}) and (\ref{e132}):
$$
\frac{x^t Qx}{\sum_i x_i^2} =
\frac{\sum_{(i,j) \in \overline{E}}  (x_i-x_j)^2} {\sum_i x_i^2}
=\frac{ \sum_{(i,j) \in \overline{E}}  (x_i-x_j)^2 
\sum_{(i,j) \in \overline{E}} h^2(i,j) (x_i +x_j)^2}
{\sum_i x_i^2 \sum_{(i,j) \in \overline{E}} h^2(i,j) (x_i +x_j)^2}
$$
$$
\geq \frac{(\sum_{(i,j) \in \overline{E}} h(i,j) (x_i^2-x_j^2))^2}
{(4+2c^2) (\sum_i x_i^2)^2}
\geq \frac{c^2}{4+2c^2},
$$
completing the proof.  $\Box$

Returning to the graph $G_1=(C_m^d)_1$  defined in the introduction, note that
it is
a $2d$-regular graph on $m^d$ vertices. If the adjacency
matrix of $C_m$ is  $A$, as before, then the one of $G_1$ is the sum of
$d$ terms, each of which is a tensor product of $d-1$ copies of
$I$ and one copy of $A$. Thus, here, too, we can use the vector
$x^{\otimes d}$, where $x$ is as in Lemma \ref{l4}, and prove the following.
\begin{lemma}
\label{lv1}
There exists a set of vertices $W$ of $G_1$ that contains no nontrivial
cycles so that its vertex boundary in $G_1$ is of size $c|W|$, and
$\frac{c}{1+c} \leq 4 \sqrt d \sin (\frac{\pi}{2m}) \leq 2 \pi \frac{\sqrt d}{m}.$
\end{lemma}
{\bf Proof:}\,
The Laplace matrix $L$ of $G_1$ and the vector $x^{\otimes d}$ above satisfy
$$
\frac{(x^{\otimes d})^t L x^{\otimes d}}{||x^{\otimes d}||^2 }
=2d-2d \cos (\pi/m)=4d \sin^2 \frac{\pi}{2m}.
$$
By Theorem \ref{t13} this implies that there is a set of vertices $W$ containing
no vertices with any coordinate being $m$ so that if $|N(W)-W|=c |W|$, then
$\frac{c^2}{4+2c^2} \leq 4d \sin^2 \frac{\pi}{2m}.$ The desired result follows, as
$$
\frac{1}{4}(\frac{c}{1+c})^2 \leq \frac{c^2}{4+2c^2} \leq 4d \sin^2 \frac{\pi}{2m}
< \frac{d \pi^2}{m^2},
$$
completing the proof.
$\Box$
\vspace{0.3cm}

\noindent
{\bf Proof of Theorem \ref{t11}:}\, 
The proof is very similar to that of Theorem \ref{t1}. 
Let $W$ be as in Lemma \ref{lv1}, 
let $v_1, v_2, \ldots$ be a random sequence
of vectors in $Z_m^d$, and define $W_i=v_i+W=\{v_i+w: w \in W\}$ where
addition is taken modulo $m$ in each coordinate. 
By symmetry, the induced subgraph of
$G_1$ on each $W_i$ is isomorphic to that on $W$ and hence contains 
no nontrivial cycle. Obviously, with probability $1$ there exists a finite
$s$ so that the union $\cup_{i=1}^s (W_i \cup N(W_i))$ covers all vertices of $G_1$. 
For each $ i $, define $B_i=(N(W_i)-W_i) -\cup_{j < i} (W_j \cup N(W_j))$.
The union of all the sets $B_i$ is 
a vertex spine, as each cycle that uses no vertex of this union is
contained in 
a single set $W_i$. We claim that for each fixed vertex $v$
of  $G_1$, the probability that $v$ belongs to the above union is 
$\frac{c}{1+c}$. Indeed, if
$i$ is the smallest $j$ so that $v \in (W_j \cup N(W_j))$, 
$v$ is a uniform random vertex of $W_i \cup N(W_i)$ and the probability that it lies in
$N(W_i)-W_i$ is thus precisely $c/(1+c)$, as claimed. By linearity of expectation, the
expected size of the union of all sets $B_i$ is $\frac{c}{1+c} m^d$, and the 
desired result follows.
$\Box$

\noindent
{\bf Remarks:}\,
\begin{itemize}
\item
A simple computation shows that for large $m$ and $d$, the expression
$\frac{2 \mu}{3^d-1}$ in Theorem \ref{t1} is  at most $(1+o(1)) \sqrt
{\frac{8 }{3}} \frac{\pi \sqrt {d}}{m}.$
\item
By the remark following the proof of Theorem \ref{t3}, the set $W$ 
in  Lemma \ref{l4}  is a level set of the vector $x^{\otimes d}$
(or equivalently, the vector obtained from it by squaring each
coordinate.)
\item
Theorem \ref{t1}  gives an alternative proof of the main 
result of Raz \cite{Ra}, showing that in the parallel repetition
for the maxcut  game on an odd cycle, $\Theta(m^2)$ repetitions
are required to ensure a value smaller than $1/2$. 
\item
The assertion of Theorems \ref{t1} and \ref{t11} can be extended to powers of other
Cayley graphs. 
\end{itemize}

\section{The continuous case}
The discrete results above  have a continuous analogue.
Let $d$ be a dimension, and let 
$\TT^d = \left. \RR^d \right / \ZZ^d$ be the unit torus.
The  torus $\TT^d$, which will be identified as a 
set with $[0,1)^d$, inherits the Riemannian structure 
of $\RR^d$. Recall that we write $\mbox{Vol}_d$ and $\mbox{Vol}_{d-1}$ for the
$d$-dimensional and $(d-1)$-dimensional Hausdorff measures 
on the unit torus $\TT^d$. A {\it body} here means a 
non-empty compact set 
that equals the closure of its interior.
A {\it smooth} function or surface always means here $C^{\infty}$-smooth.

\begin{lemma} 
 \label{lem_1124}
There exists a body $D \subset (0,1)^d \subset \TT^d$
with a smooth boundary, such that 
$$ \mbox{Vol}_{d-1}(\partial D) \leq 2 \pi \sqrt{d} \mbox{Vol}_d(D).$$
 \end{lemma}

\noindent {\bf Proof:} Denote 
\begin{equation}
 h = \inf_{A \subset (0,1)^d} \frac{\mbox{Vol}_{d-1}(\partial A)}
{\mbox{Vol}_d(A)}
\label{eq_1117}
\end{equation}
where the infimum runs over all bodies $A$ with a smooth boundary in $(0,1)^d$. Cheeger's inequality with Dirichlet boundary conditions on the cube states 
that for any smooth function $\vphi: [0,1]^d \rightarrow \RR$ 
that vanishes on the boundary,
\begin{equation}
 -\int_{[0,1]^d} \vphi \triangle \vphi \geq \frac{h^2}{4} \int_{[0,1]^d} \vphi^2,
\label{eq_1118}
\end{equation}
where $\triangle$ is the Laplacian.
For a short proof of Cheeger's inequality, see, e.g., \cite{SY}, Chapter III.
The best function to substitute in (\ref{eq_1118}) is the Laplacian 
eigenfunction $\vphi(x) = \prod_{i=1}^d \sin(\pi x_i)$, which satisfies
$\triangle \vphi = -d \pi^2 \vphi$. From (\ref{eq_1118}) we thus learn that
$d \pi^2 \leq h^2 / 4$, and the lemma follows. $\square$ 

\noindent {\bf Remark.} The set $D \subset (0,1)^d$ in Lemma \ref{lem_1124}
may be chosen to be convex. In fact, as the proof of Cheeger's inequality 
 shows, the set $D$ may be chosen to be a level set of the 
concave function $\log \vphi(x) =  \sum_{i=1}^d \log \sin(\pi x_i)$.

\medskip

\noindent {\bf Proof of Theorem \ref{thm_1053}:} 
Let $v_1,v_2,\ldots \in \TT^d$ be a sequence of independent 
random vectors, uniformly distributed in the torus $\TT^d$. 
Let $D$ stand for the body from Lemma \ref{lem_1124}, and write $D_i = v_i + D$, where addition is carried in the group $\TT^d$.
Consider the disjoint union
$$ S = \bigcup_{i=1}^{\infty} S_i, \ \ \ \ \text{where} \ \ \ \ \ S_i =  
\partial D_i \setminus \cup_{j=1}^{i-1} D_j. $$
Since $D$ has a non-empty interior, then with probability one,
$\TT^d$ is the union of finitely many $D_j$'s. Thus, if a loop in $\TT^d$
does not intersect $S$, it must be contained in $D_i \subset v_i + (0,1)^d$
for some $i$, and hence it is contractible. Consequently, $S$ is a compact spine with probability one. It remains to show that $\EE \mbox{Vol}_{d-1}(S) \leq 2 \pi \sqrt{d}$. 
Note that $S$ is a finite union of $S_i$'s, and each $S_i$
is a relatively open subset of the smooth hypersurface $\partial D_i$.
Therefore,
\begin{equation}
 \mbox{Vol}_{d-1}(S) = 
\lim_{\eps \rightarrow 0^+}
(2 \eps)^{-1} \mbox{Vol}_d \left( \bigcup_{i=1}^{\infty} (S_i)_{\eps} \right) =
\lim_{\eps \rightarrow 0^+}
\eps^{-1} \mbox{Vol}_d 
\left( \bigcup_{i=1}^{\infty} \left( (S_i)_{\eps} \cap D_i \right) \right),
\label{eq_329}
\end{equation}
where $(S_i)_{\eps}$ is the set of all points in $\TT^d$ whose geodesic 
distance from $S_i$ is smaller than $\eps$. 
Fix a point $x \in \TT^d$. Then 
\begin{equation}
 \lim_{\eps \rightarrow 0^+}
\EE \, \eps^{-1} \mbox{Vol}_d \left( \bigcup_{i=1}^{\infty} \left( (S_i)_{\eps} \cap D_i \right) \right)
=  \lim_{\eps \rightarrow 0^+}   \eps^{-1} \bP \left( x \in \bigcup_{i=1}^{\infty} \left( (S_i)_{\eps} \cap D_i \right) \right) \label{eq_326}
\end{equation}
 There exists a minimal index $i$ such that $x \in D_i$.
Let $\ell$ be this minimal index (so $\ell$ is a random variable). 
The crucial observation is that $x - v_{\ell}$ is distributed uniformly in $D$.
Hence,  we may continue (\ref{eq_326}) with 
\begin{equation} 
\label{eq_330}
= 
\lim_{\eps \rightarrow 0^+}   \eps^{-1} \bP \left( x \in \bigcup_{i=\ell}^{\infty} \left( (S_i)_{\eps} \cap D_i \right) \right)
\leq \lim_{\eps \rightarrow 0^+}  \eps^{-1} \bP \left( x \in (\partial D_\ell)_{\eps} \right) 
= \frac{\mbox{Vol}_{d-1}(\partial D_\ell)}{\mbox{Vol}_d(D_\ell)} \leq 2 \pi \sqrt{d}, 
\end{equation}
according to Lemma \ref{lem_1124}, since $\bP( x \in (\partial D_\ell)_{\eps})
= \bP(x - v_{\ell} \in (\partial D)_{\eps})$.
 From (\ref{eq_329}) and Fatou's lemma,
$$ \EE \, \mbox{Vol}_{d-1}(S) \leq \lim_{\eps \rightarrow 0^+}
\EE \, \eps^{-1} \mbox{Vol}_d \left( \bigcup_{i=1}^{\infty} \left( (S_i)_{\eps} \cap D_i \right) \right) \leq 2 \pi \sqrt{d}, $$
where the last inequality follows from (\ref{eq_326}) and (\ref{eq_330}). The proof is
complete. %Hence there exists a spine $S$ with $\mbox{Vol}_{d-1}(S) \leq 2 \pi \sqrt{d}$. 
%The compactness of $S$ is obvious.
$\Box$
%\newpage

\noindent {\bf Remarks:} 
\begin{itemize}

\item A spine $S \subset \TT^d$ is called regular if it is contained
in a finite union of smooth hypersurfaces in $\TT^d$.
A spine $S \subset \TT^d$ 
is minimal if for any $x \in S$
and $\eps > 0$, the set $S \setminus B(x, \eps)$ is no longer a spine,
where $B(x,\eps)$ is the open ball of radius $\eps$ about $x$.
By Zorn's lemma, for any compact spine $S$ there exists 
a minimal sub-spine $S^{\prime} \subset S$. We may thus assume that the 
spine $S$ in Theorem \ref{thm_1053} is minimal and regular. 
When $S$ is a minimal regular spine, the set $\TT^d \setminus S$
is necessarily connected, and since it intersects all non-contractible
loops, it is simply-connected. Note that for a minimal spine $S$,
the set $\bar{S} = \{ x \in \RR^d ; x \mod \ZZ^d \in S \}$
is the boundary of a $\ZZ^d$-periodic tiling of $\RR^d$ with connected cells of volume one.

\item Theorem \ref{thm_1053} is tight, up to the value of the constant $2 \pi$. Indeed, 
suppose $S \subset \TT^d$ is a minimal spine. Consider 
the set $\bar{S} = \{ x \in \RR^d ; x \mod \ZZ^d \in S \}$, and pick a connected component $C$ of $\RR^d \setminus \bar{S}$. Then $\mbox{Vol}_d(C) = 1$ 
and $\mbox{Vol}_{d-1}(\partial C) = 2 \mbox{Vol}_{d-1}(S)$. 
By the classical isoperimetric inequality in $\RR^d$,
$$ \mbox{Vol}_{d-1}(S) = \mbox{Vol}_{d-1}(\partial C) / 2 \geq \kappa_d 
\mbox{Vol}_d(C)^{(d-1)/d} / 2 = \kappa_d / 2 $$
where $\kappa_d = d \sqrt{\pi} \Gamma(1 + d/2)^{-1/d} \geq \sqrt{d} \left( \sqrt{2 \pi e} + o(1) \right)$.

 \item Our proof uses very few properties of the torus. 
A straightforward generalization of Theorem \ref{thm_1053} 
might read as follows:
Suppose a Lie group $G$ acts transitively by isometries on a 
simply-connected Riemannian manifold  
$\Omega$ (in our case $G = \Omega = \RR^d$).
Let $\Gamma$ be a discrete, co-compact subgroup (in our case $\Gamma = \ZZ^d$), and let $T \subset \Omega$
be a fundamental domain (in our case $T = [0,1)^d$). Assume that $T$ is simply connected,
and write $\lambda$ for the minimal eigenvalue of 
minus the Laplacian with Dirichlet boundary conditions on $T$. 
Then, 
there exists a compact spine in $\Omega / \Gamma$ whose surface area is at most $2 \sqrt{\lambda}$. 

Note that there clearly exists a trivial spine in  $\Omega / \Gamma$
whose area is at most $\mbox{Vol}_{d-1}(\partial T)$. Only in the case where
$\sqrt{\lambda} << \mbox{Vol}_{d-1}(\partial T)$ we 
obtain a non-trivial conclusion. 

\item A short argument leading from the continuous 
Theorem \ref{thm_1053} to the discrete Theorem \ref{t1} (with a slightly
worse constant) appears in \cite{FKO}, Theorem 3.1.

\end{itemize}

\noindent
{\bf Acknowledgment:} We thank Guy Kindler, Anup Rao and Avi Wigderson for helpful comments.

\end{document}